\documentclass{amsart}
\usepackage{amsmath,amsthm}
\usepackage{amssymb}
\usepackage{euscript}
\usepackage{latexsym}
\usepackage[all]{xy}
\usepackage{verbatim}

\newcommand{\THH}{\operatorname{THH}}

\newcommand{\hocolim}{\operatornamewithlimits{hocolim}}
\newcommand{\colim}{\operatornamewithlimits{colim}}

\newcommand{\cy}{{\text{\rm cy}}}

\newcommand{\id}{\operatorname{id}}
\newcommand{\pr}{\operatorname{pr}}

\newcommand{\xto}{\xrightarrow}

\begin{document}

\def\ra{\longrightarrow}
\def\la{\longleftarrow}
\newtheorem{thm}{Theorem}[section]
\newtheorem{defn}[thm]{Definition}
\newtheorem{prop}[thm]{Proposition}
\newtheorem{claim}[thm]{Claim}
\newtheorem{cor}[thm]{Corollary}
\newtheorem{lem}[thm]{Lemma}
\newtheorem{rem}[thm]{Remark}

\title{Topological Hochschild Homology of Twisted Group Algebras}
\author{Daniel J. Vera}
\date{May 31, 2006}

\maketitle

\begin{abstract}
We show the topological Hochschild homology spectrum of a twisted group 
algebra $\THH(A^{\tau}[G])$ is the Thom spectrum associated to a 
parametrized orthogonal spectrum $E(A,G)$.  We then analyze the structure of the parametrized orthogonal spectrum $E(A,G)$ and show that it is locally trivial.
\end{abstract}

\section{Introduction}

Let $G$ be a group and $A$ be a commutative ring.  The topological Hochschild homology spectrum of the group algebra $A[G]$ 
is determined by the cyclic bar construction of the group $G$ and the topological Hochschild spectrum of the ring $A.$  
More precisely, Hesselholt and Madsen proved \cite[Thm.~7.1]{hm} there is a stable equivalence of orthogonal spectra $$\THH(A)\wedge N^{\cy}(G)_{+}\stackrel{\sim}{\ra} \THH(A[G]).$$
This paper generalizes this result to a {\it twisted group algebra} $A^{\tau}[G]$.  Let $G$ act on $A$ from the left through ring automorphisms.  Then as an $A$-module, $A^{\tau}[G]=A[G],$ but the multiplication
is given by $ag \cdot a' g' = ag(a') \cdot gg'.$  The result is stated in terms of a variant $\THH^{g}(A)$ of the topological Hochschild spectrum that we
describe in \S\ref{defofE(A,G)}.  We then prove the following result.
\begin{thm}\label{mainresult} Let $G$ be a group that acts on a ring $A,$ and let $A^{\tau}[G]$ be the twisted group algebra.  Then there is a stable equivalence of orthogonal spectra  $$\Phi : \bigvee\limits_{\langle g \rangle} EG_{+} \wedge_{C_{G}(g)} \THH^{g}(A)\stackrel{\sim}{\ra} \THH(A^{\tau}[G]),$$ where 
the wedge-sum on the left hand side ranges over the conjugacy classes of elements of $G$.  The map $\Phi$ depends on a choice of representative $g \in \langle g \rangle$ of every conjugacy class of elements in $G.$
\end{thm}

We first prove a form of this equivalence that is independent of the choice of conjugacy class representatives.  Let $\lambda$ be a finite dimensional inner product space and let $S^\lambda$ be the one point compactification of $\lambda$.  The topological Hochschild spectrum is the orthogonal spectrum defined as the geometric realization $$\THH(A)_\lambda = \THH(A;S^{\lambda}) = \big| [k] \mapsto \THH(A;S^{\lambda})[k] \big|$$ of a cyclic orthogonal spectrum $\THH(A;S^{\lambda})[k]$.  For each non-negative integer $k,$ we define a parametrized orthogonal spectrum $E(A,G)[k]$ over the cyclic bar construction $N^{\cy}(G)[k]$ with $\lambda$-th space given by $$E(A,G)[k]_{\lambda} = E(A,G;S^{\lambda})[k] = \THH(A;S^{\lambda})[k]\times N^{\cy}(G)[k].$$
We then define cyclic operators $d^{\tau}_{i,E}, s^{\tau}_{i,E},$ and $t^{\tau}_{k,E}.$  For $0\leq i\leq k,$ the degeneracy operator $s^{\tau}_{i,E}$ is given as the product of the $i$-th degeneracy operator of the cyclic pointed space $\THH(A,S^{\lambda})[-]$ and the $i$-th degneracy operator of the cyclic set $N^{\cy}(G)[-],$ respectively.  Similarily, the cyclic operator $t^{\tau}_{k,E}$ is the product of the cyclic operator of $\THH(A;S^{\lambda})[-]$ and the cyclic operator of $N^{\cy}(G)[-].$  The face operators, however, are replaced by {\it twisted} face operators defined in \S\ref{defofE(A,G)}.  At each simplicial level $k,$ $E(A,G)[k]$ is a parametrized orthogonal spectrum over 
$N^{\cy}(G)[k]$ but it is not a cyclic object in a category of parametrized orthogonal spectra over a fixed base. 
Fixing $\lambda$ and letting $k$ vary, we have a cyclic space, $E(A,G)[-]_{\lambda}.$  The geometric realization of $E(A,G)[-]_{\lambda}$ is 
$$E(A,G)_{\lambda} = \big| [k]\mapsto E(A,G)[k]_{\lambda}\big|.$$   At each level $k,$ the projection map of the parametrized space $E(A,G)[-]_{\lambda}$ onto $N^{\cy}(G)[-]$ commutes with the operators $d^{\tau}_{i,E},$ $s^{\tau}_{i,E},$ and $t^{\tau}_{i,E}.$  Thus for varying $k$ the projection maps form a cyclic map.  Similarly the level $k$ section maps form a cyclic map.  Since the geometric realization is a functor, it follows that $E(A,G)_{\lambda}$ is a parametrized space over $N^{\cy}(G).$  Letting $\lambda$ vary then gives a parametrized orthogonal spectrum $E(A,G)$ over $N^{\cy}(G)$ in the sense 
of Definition~\ref{paramspecdef}.  
The space $E(A,G)_{\lambda}$ is also a bundle over $N^{\cy}(G)$ with the fiber over a vertex $g \in \langle g \rangle$ denoted by $E(A,G)^{g}_{\lambda}.$  The fiber $E(A,G)^{g}_{\lambda}$ is itself obtained as the realization of a cyclic pointed space $E(A,G)^{g}[-]_\lambda.$  We discuss this in \S\ref{THHcyc}. 

In general, if $\mathcal{T}_{\ast}$ denotes the category of pointed spaces and $\mathcal{T}_{B}$ the category of parametrized spaces over a base $B$, there exists
an adjoint pair of functors $(f_{!}, f^{\ast})$ between parametrized spaces and pointed spaces, $f_{!}:\mathcal{T}_{B} \ra \mathcal{T}_{\ast}$, the associated Thom space, and $f^{\ast}:\mathcal{T}_{\ast}\ra \mathcal{T}_{B},$ the change of base functor.  May and Sigurdsson show that applying the Thom space functor $f_{!}$ levelwise to a parametrized orthogonal spectrum gives us an orthogonal spectrum  \cite[Thm. ~11.4..1]{maysigurdsson}.  Since the functor $f_{!}$ is a left adjoint, it preserves colimits.  Thus for the parametrized spectrum $E(A,G),$ the Thom spectrum $f_{!}E(A,G)$ is given by the realization of a cyclic orthogonal spectrum whose orthogonal spectrum in 
simplicial degree $k$ has $\lambda$-th space given by $\THH(A;S^{\lambda})[k]\wedge N^{\cy}(G)[k]_{+}.$
\begin{thm}\label{secondmainresult}  There exists a canonical stable equivalence of orthogonal spectra $$\Psi : f_{!}E(A,G) \stackrel{\sim}{\ra} \THH(A^{\tau}[G]).$$ 
\end{thm}  

By Connes' theory of cyclic sets, the realization of a cyclic set is a ${\mathbb T}$-space.  Since the topological Hochschild homology spectrum is defined as the realization of a cyclic orthogonal spectrum, it is equipped with an action of the circle; see Loday's book \cite{loday}.  In particular, the spectrum $\THH(A^{\tau}[G])$ is an orthogonal $\mathbb{T}$-spectrum as defined by Mandell and May in \cite{mandellmay} and Theorem~\ref{secondmainresult} can be extended to an equivalence of orthogonal $\mathbb{T}$-spectra.  However, the fixed
points of an equivariant Thom spectrum such as $f_! E(A,G)$ are
notoriously difficult to understand; compare the calculation by Araki \cite{araki} of the $\mathbb{Z}/2\mathbb{Z}$-fixed set $MR^{\mathbb{Z}/2\mathbb{Z}}$ of the Thom $\mathbb{Z}/2\mathbb{Z}$-spectrum that corresponds to Real $K$-theory.
  
The proof of Theorem~\ref{secondmainresult} employs the definition of $\THH$ and the base change functors $(f_{!},f^{\ast})$ from parametrized homotopy theory.  It is independent of the twisting in the simplicial structure in that the equivalence does not depend on a preselected choice of element $g \in \langle g \rangle$ as is the case for Theorem ~\ref{mainresult}.  The proof of Theorem ~\ref{mainresult} involves an explicit analysis of the cyclic structure of the cyclic bar construction and is inspired by a study of the linear case of the Hochschild homology.  Another formula for the (ordinary) Hochschild homolgy of twisted group algebras is given by Fe\u{\i}gin and Tsygan in \cite[\S 4]{feigintsygan}, but note that we use a different system of coordinates for the cyclic bar construction (see \S\ref{defofE(A,G)}). 

We use the following notational conventions throughout this exposition.  By a space we mean a compactly generated space (weak Hausdorff $k$-space) and by a pointed space we mean a compactly generated space with a choice of base-point.  Let $\mathcal{S}$ be the category of sets.  We denote both a given cyclic set $$X: \Lambda^{\operatorname{op}}\ra\mathcal{S}$$ and a given simplicial set $$X: \Delta^{\operatorname{op}}\ra \mathcal{S}$$ by $X[-].$  We will always state whether we are considering a simplicial set or a cyclic set so that no confusion arises.  We also let $X$ denote both the geometric realization of a simplicial set $X[-]$ and the geometric realization of the cyclic set $X[-],$ defined as the geometric realization of the underlying simplicial 
set.  Finally, whenever we make use of the functors $f_{!}$ and $f^{\ast}$, $f$ will always be the map sending the base to the one point space unless explicitly stated otherwise.

It is a great pleasure for the author to thank his thesis advisor Lars Hesselholt for his invaluable guidance in this research and Haynes Miller for all of his helpful comments about this paper.  This research was completed at the Massachusetts Institute of Technology and Nagoya University.  The author also thanks the MIT Department of Mathematics for their support during his graduate studies, and thanks the Graduate School of Mathematics at Nagoya University for their kind hospitality and support through the COE program at Nagoya University during his visit in the Spring of 2005. 

\section{Parametrized Spaces}\label{paramspaces}
 
We discuss parametrized spaces and define the adjoint functors $f_{!}$ and $f^{\ast}$ that we use in Theorem ~\ref{mainresult} and Theorem ~\ref{secondmainresult}.  We will also define a right adjoint to $f^{\ast},$ denoted $f_{\ast}.$ 
The reference for this material is \cite{maysigurdsson}.  \begin{defn}  Let $B$ be a fixed base space.  A parametrized space $X$ over $B$ consists of a space $X$ together with projection and section maps, $p:X\ra B$ and $s:B\ra X$ respectively, such that $p\circ s = \id_{B}.$
\end{defn}  The category of parametrized spaces $\mathcal{T}_{B}$ over a fixed base space $B$ has objects parametrized spaces over $B.$  The morphisms of $\mathcal{T}_{B}$ are maps of the total spaces that commute with both the section and projection maps, or commutative diagrams:
\[\xymatrix{&X\ar[dd]^{f}\ar[dr]^{p_{1}}\\
B\ar[ur]^{s_{1}}\ar[dr]_{s_{2}}& &B \\
&Y\ar[ur]_{p_{2}}}\]
We construct the {\it parametrized mapping space} $F_{B}(-,-)$ and {\it parametrized smash product} to make $\mathcal{T}_{B}$ a closed symmetric monoidal category with unit $S_{B}^{0} = B\times S^{0}.$  The zero object $\ast_{B}$ in $\mathcal{T}_{B}$ is the space $B$ with projection and section maps given by the identity.  Let $\mathcal{T}_{B}(X,Y)$ be the set of all morphisms in $\mathcal{T}_{B}$ from $X$ to $Y.$  We topologize this set as a subspace of the space of all unbased maps of unbased total spaces $X\ra Y.$  We note that the space $\mathcal{T}_{B}(X,Y)$ is a based
space with basepoint the map $s_{2}\circ p_{1}:X\ra Y.$  This is the unique map factoring through $\ast_{B}$ in $\mathcal{T}_{B}.$  Thus the category $\mathcal{T}_{B}$ is enriched over $\mathcal{T}_{\ast}.$  It is also {\it based topologically bicomplete} in the sense of \cite[~\S 1.2]{maysigurdsson}.  Constructing $F_{B}(-,-)$ requires us to first construct the
{\it unbased} parametrized mapping space $Map_{B}(-,-)$.  To do this, we introduce a subtle preliminary notion.  For a space $Y \in \mathcal{T}$ (unbased), the {\it partial map classifier} is $\tilde{Y} = Y \cup \{\omega\}$ where $\omega$ is a disjoint basepoint.  It is topologized as the space with basis $\{ U \cup \{\omega\} : U \in \mathcal{U} \}$ where $\mathcal{U}$ is a basis for $Y$ \cite[Def. ~1.3.10]{maysigurdsson}.  We note that the point $\omega$ is not closed and $\tilde{Y}$ is not weak Hausdorff \cite[Def. ~1.1.1]{maysigurdsson}.  Also the closure of $\{ \omega \}$ is all of $\tilde{Y}.$  The point $\omega$ is analagous
to a generic point of a variety.  The space $\tilde{Y}$ is known as the {\it partial map classifier} because of the bijective correspondence between maps $f:A \ra Y$ with $A\subseteq X$ a closed subset, and corresponding maps $\tilde{f}: X\ra \tilde{Y}$ defined by $$\tilde{f}(x) = \begin{cases} \omega & \text{if $x\notin A$}\\
f(a) & \text{if $x\in A.$} \end{cases}$$ 
For a space $p:X\ra B$ over $B$ we define a map $\xi: B\ra Map(X, \tilde{B})$ by $$\xi(b) (x) = \begin{cases} b  & \text{if $x\in X_{b}$}\\
\omega & \text{otherwise.} \end{cases}$$  Here $X_{b}$ denotes the fiber $p^{-1}(b)$ over $b\in B.$ The map $\xi$ is the adjoint of the map $\tilde{f}: X\times B \ra B$ obtained as $f:\Delta^{-1}(p\times \id (X\times B))\ra B$ where $p\times \id:X\times B \ra B\times B$ and $\Delta:B\ra B\times B$ denotes the diagonal embedding.  We now have the following.
\begin{defn} \cite[Def. ~1.3.11]{maysigurdsson} Let $p:X\ra B$ and $q:Y\ra B$ be spaces over $B$ and $Map(X,Y)$ the space of unbased maps from $X$ to $Y$.  Then $Map_{B}(X,Y)$ is defined to be the pullback of the following diagram,
\[\xymatrix{Map_{B}(X,Y)\ar[d]\ar[r] & Map(X, \tilde{Y})\ar[d]^{\tilde{q}_{\ast}} \\
B \ar[r]^(.4){\xi} & Map(X, \tilde{B}).
}\] 
\end{defn}  We note that as a point-set, $$Map_{B}(X,Y) = \coprod_{b\in B} Map(X_{b}, Y_{b}).$$  We may now define the {\it parametrized mapping space},
\begin{defn}  \cite[Def. ~1.3.16]{maysigurdsson} The parametrized mapping space $F_{B}(X,Y)$ of two parametrized spaces $X$ and $Y$ is the pullback of the following diagram,
\[\xymatrix{F_{B}(X,Y)\ar[d]\ar[rr] && Map_{B}(X,Y) \ar[d]^{(s_{1})_{\ast}}\\
B\ar[r]^(.475){s_{2}}& Y\ar[r]^(.275){\sim}& Map_{B}(B,Y)
}\] where $s_{1}$ and $s_{2}$ are the sections of $X$ and $Y$, respectively, and $Y \stackrel{\sim}{\ra} {Map}_B(B,Y)$ is the canonical isomorphism. 
\end{defn}
Again we note that as a point-set, $$F_{B}(X,Y) = \coprod_{b\in B}F(X_{b},Y_{b}).$$  The parametrized mapping space is thus the subspace of $Map_{B}(X,Y)$ consisting of maps that restrict to based maps between the fibers $X_{b}$ and $Y_{b}$ with respective basepoints $s_{1}(b)$ and $s_{2}(b).$

Given a map of spaces $f:A \ra B$ we define a pair of adjoint functors, $$f_{!}:\mathcal{T}_{A} \ra \mathcal{T}_{B}$$ and $$f^{\ast}:\mathcal{T}_{B} \ra \mathcal{T}_{A}$$ by the following pushout and pullback diagrams, respectively.
\[\xymatrix{A \ar[d]_{s} \ar[r]^{f} & B \ar[d] && f^{\ast}Y \ar[d] \ar[r] & Y \ar[d]^{p}\\
X \ar[r] & f_{!}X && A \ar[r]^{f} & B }\] where $X\in \mathcal{T}_{A}$ and $Y\in \mathcal{T}_{B}$.  
Of particular interest is the example where $f:B\ra \ast$ is the map sending $B$ to the one-point space.
In this case, for a parametrized space $X$ over $B$, one obtains the pointed space $f_{!}(X) = X/s(B)$ with basepoint provided by the class of the section.  Similarly, given a pointed space $Z$ with basepoint $z_{0}\in Z$,
one has the parametrized space $f^{\ast}(Z)$ over $B$ with total space $Z \times B$ and projection provided by projecting onto the second factor.  The section $s:B\ra Z\times B$ is then defined by
$s(b)=(z_{0},b).$  The functor $f^{\ast}$, has a right adjoint, $f_{\ast}:\mathcal{T}_{A} \ra \mathcal{T}_{B}$ defined as follows.  Let $\iota : B\ra Map_{B}(A,A)$ be the adjoint of the map $A\times_{B}B \ra A$ sending $(a, f(a))\mapsto a.$  Then for $X\in \mathcal{T}_{A}$ we define $f_{\ast}X$ as the pullback
\[\xymatrix{f_{\ast}X \ar[d]\ar[r]&Map_{B}(A,X)\ar[d]\\
B\ar[r]^(.3){\iota}& Map_{B}(A,A).
}\]  In the case where $f:B\ra \ast, f_{\ast}X$ is the space of all sections on $X$ with basepoint the section $s:B\ra X.$    

Within the category of parametrized spaces, we can define parametrized versions of the wedge and smash products by taking the usual wedge and smash products fiberwise.   More precisely, we have the following definition.  
\begin{defn}  \cite[Def. ~1.3.8]{maysigurdsson} Let $X, Y \in \mathcal{T}_{B}.$
\begin{itemize}
\item[(1)] The product of spaces $X$ and $Y$ over $B$, $X\times_{B}Y$, is the pullback of the following diagram,  
\[\xymatrix{ X \times_{B}Y \ar[d]\ar[r]&X \ar[d]^{p_{1}}\\
Y\ar[r]_{p_{2}} & B }\]

\item[(2)] The wedge of spaces $X$ and $Y$ over $B$, $X \vee_{B} Y$, is the pushout obtained from the following diagram, 
\[\xymatrix{  B \ar[d]_{s_{2}} \ar[r]^{s_{1}}&X \ar[d]\\
Y\ar[r] & X \vee_{B}Y }\]

\item[(3)] Finally, we have the inclusion map $X \vee_{B}Y \ra X\times_{B}Y$ defined by sending $x \mapsto (x, s_{2}p_{1}(x))$ and $y \mapsto (s_{1}p_{2}(y),y)$.  It is easy to check that it is well-defined.  We then define the smash product of $X$ and $Y$ over $B$, $X \wedge_{B} Y$, to be the pushout,
\[\xymatrix{X \vee_{B}Y \ar[d] \ar[r]& X \times_{B}Y \ar[d]\\
B \ar[r] & X \wedge_{B}Y. }\]
\end{itemize}
\end{defn}  Since every fiber has a basepoint given by the section, each of the above constructions for the product, wedge product, and smash product, is given the parametrized space structure that gives us fiberwise product, wedge product, and smash product, respectively.  

\section{Parametrized Orthogonal Spectra}\label{paramspectra}

We will now define parametrized orthogonal spectra and discuss how the base change functors of \S\ref{paramspaces} extend to parametrized orthogonal spectra.  Furthermore, if we are given a simplicial (or cyclic) space $B[-]$ and for each non-negative integer $k$ a parametrized orthogonal spectrum $E[k]$ over $B[k]$, we will construct a parametrized orthogonal spectrum $E$ over the geometric realization $B$ of $B[-]$.  

Let $\lambda$ be a finite dimensional real inner product space and let $S^{\lambda}$ be the one-point compactification of $\lambda.$  A {\it topological category} is a category enriched in the symmetric monoidal category of pointed spaces and smash product.  Let $\mathcal{I}$ be the topological category with objects all finite dimensional real inner product spaces $\lambda$ and morphism spaces given by the pointed space of linear isometries from 
$\lambda$ to $\lambda ',$ $$\operatorname{Hom}_{\mathcal{I}}(\lambda, \lambda') = \mathcal{O}(\lambda, \lambda ')_+.$$  Let 
\[\xymatrix{E(\lambda, \lambda ') \ar[d]\ar@{^{(}->}[r]& \mathcal{O}(\lambda, \lambda ') \times \lambda' \ar[d]^{pr_{1}}\\
\mathcal{O}(\lambda, \lambda') \ar@{=}[r]& \mathcal{O}(\lambda, \lambda ')
}\] be the sub-bundle of pairs $(f,x)$ such that $x\in \lambda ' - f(\lambda),$ the orthogonal complement.  Let $\mathcal{J}$ be the topological category with the same objects as $\mathcal{I}$ but with morphism spaces $\operatorname{Hom}_{\mathcal{J}}(\lambda, \lambda')$ defined to be the Thom space of the vector bundle $E(\lambda, \lambda')$ over $\mathcal{O}(\lambda, \lambda')$.  Composition is defined $$\operatorname{Hom}_{\mathcal{J}}(\lambda ', \lambda '') \wedge \operatorname{Hom}_{\mathcal{J}}(\lambda, \lambda ') \ra \operatorname{Hom}_{\mathcal{I}}(\lambda, \lambda '')$$ via
$((g,y); (f,x)) \mapsto (g\circ f, g(x) + y).$  The inclusion of the zero-section in $E(\lambda, \lambda')$ induces a map of Thom spaces $\operatorname{Hom}_{\mathcal{I}}(\lambda,\lambda') \ra \operatorname{Hom}_{\mathcal{J}}(\lambda,\lambda')$ and this map is an isomorphism if the dimensions of $\lambda$ and $\lambda'$ are equal.  These maps constitute a functor $\mathcal{I}\ra \mathcal{J}$.  A {\it pointed-topological functor} is a functor enriched over pointed spaces.  By definition an {\it orthogonal spectrum} is a pointed-topological functor $$X: \mathcal{J} \ra \mathcal{T}_{\ast}.$$  The topological Hochschild homology spectrum $\THH$ defined in \S\ref{defofE(A,G)} is an example of an orthogonal spectrum.

We recall that the category $\mathcal{T}_B$ is enriched in the symmetric monoidal category of pointed spaces and smash product.  Let $S_{B}^{\lambda} = f^{\ast}(S^{\lambda}).$  
\begin{defn}\label{paramspecdef}  A parametrized orthogonal spectrum over $B$ is a pointed-topological functor $$X: \mathcal{J}\ra \mathcal{T}_B.$$
\end{defn} (Compare \cite[Def.~11.2.3]{maysigurdsson}).  This amounts to a pointed-topological functor (that we denote by the same symbol) $$X: \mathcal{I}\ra \mathcal{T}_B$$ together with continuous natural transformations $$\sigma_{\lambda, \lambda'} : X(\lambda) 
\wedge_B S_{B}^{\lambda'} \ra X(\lambda \oplus \lambda')$$ of pointed-topological functors from $\mathcal{I} \times \mathcal{I}$ to $\mathcal{T}_B$ such that $$\sigma_{\lambda,0}: X(\lambda) \wedge_B S_{B}^{0} \ra X(\lambda \oplus 0)$$ is the canonical isomorphism, and such that the diagram
\[\xymatrix{ X(\lambda)\wedge_{B} S_{B}^{\lambda'} \wedge_{B} S_{B}^{\lambda''}\ar[d]^{\sim}\ar[rr]^{\sigma_{\lambda, \lambda'} \wedge \id} && X(\lambda \oplus \lambda') \wedge_{B} S_{B}^{\lambda''}\ar[d]^{\sigma_{\lambda \oplus \lambda', \lambda''}}\\
X(\lambda) \wedge_{B} S_{B}^{\lambda' \oplus \lambda''} \ar[rr]^{\sigma_{\lambda, \lambda' \oplus \lambda''}}&& X(\lambda \oplus \lambda' \oplus \lambda'')
 }\] commutes.  Here the left-hand vertical map is the canonical isomorphism.  We note that in Definition~\ref{paramspecdef} when $B=\ast,$ we obtain the usual definition of an orthogonal spectrum.  

We recall from \S\ref{paramspaces} that a map of spaces $f:A \ra B$ gives rise to adjoint functors $f_!:\mathcal{T}_A \ra \mathcal{T}_B, f^\ast:\mathcal{T}_B \ra \mathcal{T}_A$ and $f_{\ast}:\mathcal{T}_A \ra \mathcal{T}_B,$ which is right adjoint to $f^{\ast}.$  For the category of parametrized orthogonal spectra over $A$ and the category of parametrized orthogonal spectra over $B$, levelwise application of the functors $f_!, f^{\ast},$ and $f_{\ast}$ gives rise to adjoint functors (that we also denote $f_!, f^{\ast}$ and $f_{\ast}$) between the categories of parametrized orthogonal spectra over $A$ and parametrized orthogonal spectra over $B$ \cite[Thm.~11.4.1]{maysigurdsson}.  Furthermore, if $g: B \ra C$ is another map of spaces, then there are canonical isomorphisms, $$(g\circ f)_{!}E \stackrel{\sim}{\ra} g_!f_!E,\phantom{gf} (g\circ f)^{\ast} E' \stackrel{\sim}{\ra}f^\ast g^\ast E', \phantom{gf} (g\circ f)_{\ast} E \stackrel{\sim}{\ra} g_\ast f_\ast E,$$ where $E$ is a parametrized orthogonal spectrum over $A$ and $E'$ is a parametrized orthogonal spectrum over $C$. 

Let $B[-]$ be a simplicial (or cyclic) space.  Suppose for all non-negative integers $k$, we have a parametrized orthogonal spectrum $E[k]$ over $B[k]$, and for every map $\theta : [m] \ra [n]$ in the simplicial index category, we have a map of parametrized orthogonal spectra over $B[n],$ $$\theta_E : E[n] \ra \theta_{B}^{\ast} E[m].$$  Here $\theta_{B}^{\ast}$ is the base-change functor associated with the map of spaces $\theta_B : B[n] \ra B[m].$  We shall require that if $\theta: [m] \ra [n]$ and $\theta': [n] \ra [p]$ are two composable maps in the simplicial index category, then the following diagram of parametrized orthogonal spectra over $B[p]$ commutes:
\[\xymatrix{E[p]\ar[d]^{\theta'_{E}} \ar[r]^(.325){(\theta' \circ \theta)_{E}}& (\theta_B \circ \theta'_{B})^{\ast} E[m]\ar[d]^{\sim}\\
{\theta'_{B}}^{\ast} E[n]\ar[r]^(.45){{\theta'_{B}}^{\ast}\theta_{E}} & {\theta'_{B}}^{\ast}\theta_{B}^{\ast}E[m].
}\]  Here the right-hand vertical map is the canonical isomorphism.  We recall that the realization $B$ of $B[-]$ is defined to be the coequalizer
$$\xymatrix{
{ \displaystyle{\coprod_{\theta \colon [m] \to [n]}^{\phantom{\theta
      \colon [m] \to [n]}} B[n] \times \Delta^m} }
\ar[r]<.7ex>^(.55){f_B} \ar[r]<-.7ex>_(.55){g_B} &
{ \displaystyle{\coprod_{ [k] }^{\phantom{ [k] }} B[k] \times
    \Delta^k} } \ar[r]^(.66){\epsilon_B} &
{ B } \cr
}$$ where, on the second summand indexed by $\theta: [m] \ra [n],$ the map $f_B$ is the unique map that, for every map $\theta \colon [m] \ra
[n]$ in the simplicial index category, makes the following diagram
commute:
$$\xymatrix{
{ B[n] \times \Delta^m } \ar[r]^{\theta_B \times \id}
\ar[d]^(.45){\operatorname{in}_{\theta}} &
{ B[m] \times \Delta^m } \ar[d]^(.45){\operatorname{in}_{[m]}} \cr
{ \displaystyle{\coprod_{\varphi \colon [k] \to [l]} B[l] \times \Delta^k} }
\ar[r]<1ex>^(.55){f_B} &
{ \displaystyle{\coprod_{ [k] } B[k] \times \Delta^k.} } \cr
}$$
The map $g_B$ is defined similarly as the unique map that, for every map $\theta \colon [m] \ra [n]$ in the simplicial index category, makes the following diagram commute:
$$\xymatrix{
{ B[n] \times \Delta^m } \ar[r]^{\id \times \theta_{\Delta}}
\ar[d]^(.45){\operatorname{in}_{\theta}} &
{ B[n] \times \Delta^n } \ar[d]^(.45){\operatorname{in}_{[n]}} \cr
{ \displaystyle{\coprod_{\varphi \colon [k] \to [l]} B[l] \times \Delta^k} }
\ar[r]<1ex>^(.55){g_B} &
{ \displaystyle{\coprod_{ [k] } B[k] \times \Delta^k.} } \cr
}$$   Let $\epsilon'_B = \epsilon_B \circ f_B = \epsilon_B \circ g_B.$  Let $\pr_{m,n} \colon B[n] \times \Delta^m \ra B[n]$ and $\pr_k \colon B[k] \times \Delta^k \ra B[k]$ be the canonical projections.  We then have the parametrized orthogonal spectrum $\pr_{m,n}^{\ast} E[n]$ over $B[n] \times \Delta^m$ and the parametrized orthogonal spectrum $\pr_k^{\ast} E[k]$ over $B[k] \times \Delta^k$.  Let $$\coprod_{[k]} \operatorname{in}_{[k]!} \pr_{k}^{\ast} E[k]$$ denote the coproduct of the parametrized orthogonal spectra $\operatorname{in}_{[k]!}\pr_k^{\ast} E[k]$ over $\coprod_{[k]} B[k] \times \Delta^k$ and let $$\coprod_{\theta\colon [m] \to [n]} \operatorname{in}_{\theta!} \pr_{m,n}^{\ast} E[n]$$ denote the coproduct of the parametrized orthogonal spectra $\operatorname{in}_{\theta!} \pr_{m,n}^{\ast} E[n]$ over $\coprod_{\theta\colon [m]\ra[n]} B[n] \times \Delta^m.$  We then define the parametrized spectrum $E$ over $B$ to be the following coequalizer of parametrized spectra over $B$:  $$\xymatrix{
{ \displaystyle{\epsilon_{B!}' (\coprod_{\theta \colon [m] \to
      [n]}^{\phantom{\theta \colon [m] \to [n]}} \operatorname{in}_{\theta!} \pr_{m,n}^{\ast} E[n])
} } \ar[r]<.7ex>^(.55){f_E} \ar[r]<-.7ex>_(.55){g_E} &
{ \displaystyle{\epsilon_{B!} (\coprod_{ [k] }^{\phantom{ [k] }}
    \operatorname{in}_{[k]!} \pr_{k}^{\ast} E[k] ) } } \ar[r]^(.73){\epsilon_E} &
{ E. } \cr
}$$ Here, on the summand indexed by $\theta:[m]\ra [n]$, $f_E$ is defined as follows.  First, the functor $f_{B!}$ commutes
with coproducts since it has a right-adjoint functor $f_B^*$.  We define a map $$f_E' \colon \coprod_{\theta \colon [m] \to [n]}
f_{B!} \operatorname{in}_{\theta!} \pr_{m,n}^*E[n] \ra
\coprod_{[k]} \operatorname{in}_{[k]!} \pr_k^* E[k]$$
of parametrized orthogonal spectra over $\coprod_{[k]} B[k] \times
\Delta^k$ to be the unique map such that, for every map
$\theta \colon [m] \ra [n]$ in the simplicial index category, the
following diagram of parametrized orthogonal spectra over
$\coprod_{[k]} B[k] \times \Delta^k$ commutes:
$$\xymatrix{
{ f_{B!} \operatorname{in}_{\theta!} \pr_{m,n}^*E[n] }
\ar[r]^(.52){f_{E,\theta}'} \ar[d]^(.45){\operatorname{in}_{\theta}} &
{ \operatorname{in}_{[m]!} \pr_m^* E[m] }
\ar[d]^(.45){\operatorname{in}_{[m]}} \cr
{ \displaystyle{ \coprod_{\varphi \colon [k] \to [l]} f_{B!}
    \operatorname{in}_{\varphi!} \pr_{k,l}^* E[l] } }
\ar[r]<1ex>^(.63){f_E'} &
{ \displaystyle{ \coprod_{[k]} \pr_k^* E[k] }. } \cr
}$$  Here the map $f_{E,\theta}'$ is defined to be the following composite
map:
$$\begin{aligned}
{} & f_{B!} \operatorname{in}_{\theta!} \pr_{m,n}^* E[n] \stackrel{\sim}{\ra}
\operatorname{in}_{[m]!} (\theta_B \times \id)_! \pr_{m,n}^* E[n] \cr
{} & \stackrel{\sim}{\ra} \operatorname{in}_{[m]!} \pr_m^* \theta_{B!} E[n] 
\xto{ \operatorname{in}_{[m]!} \pr_m^* \theta_E^{\#} }
\operatorname{in}_{[m]!} \pr_m^* E[m]. \cr
\end{aligned}$$
The first map is the unique natural isomorphism $f_{B!}
\operatorname{in}_{\theta!} \stackrel{\sim}{\ra} \operatorname{in}_{[m]!}
(\theta_B \times \id)_!$ that exists because
$$f_B \circ \operatorname{in}_{\theta} = (\theta_B \times \id) \circ
 \operatorname{in}_{[m]} \colon B[n] \times \Delta^m
\ra \coprod_{[k]} B[k] \times \Delta^k.$$
The second map is induced from the unique natural isomorphism of
\cite[Prop.~2.2.9]{maysigurdsson}, $(\theta_B \times \id)_! \pr_{m,n}^*
 \stackrel{\sim}{\ra} \pr_m^* \theta_{B!}$, that exists because the following
 diagram is a pull-back:
$$\xymatrix{
{ B[n] \times \Delta^m } \ar[rr]^{\theta_B \times \id}
\ar[d]^{\pr_{m,n}} &&
{ B[m] \times \Delta^m } \ar[d]^{\pr_m} \cr
{ B[n] } \ar[rr]^{\theta_B} &&
{ B[m]. } \cr
}$$
Finally, the last map is induced by the map $\theta_E^{\#} \colon
\theta_{B!} E[n] \to E[m]$ that is the adjoint of the given map
$\theta_E \colon E[n] \to \theta_B^* E[m]$.  Then $f_E$ is defined to be the map
$\epsilon_{B!}f_E'$.  To define the map $g_E$, again, we define a map
$$g_E' \colon \coprod_{\theta \colon [m] \ra [n]}
g_{B!} \operatorname{in}_{\theta!} \pr_{m,n}^*E[n] \ra
\coprod_{[k]} \operatorname{in}_{[k]!} \pr_k^* E[k]$$
of parametrized orthogonal spectra over $\coprod_{[k]} B[k] \times
\Delta^k$ to be the unique map
such that, for every map $\theta \colon [m] \ra [n]$ in the simplicial
index category, the following diagram of parametrized orthogonal
spectra over $\coprod_{[k]} B[k] \times \Delta^k$ commutes:
$$\xymatrix{
{ g_{B!} \operatorname{in}_{\theta!} \pr_{m,n}^*E[n] }
\ar[r]^(.52){g_{E,\theta}'} \ar[d]^(.45){\operatorname{in}_{\theta}} &
{ \operatorname{in}_{[n]!} \pr_n^* E[n] }
\ar[d]^(.45){\operatorname{in}_{[n]}} \cr
{ \displaystyle{ \coprod_{\varphi \colon [k] \to [l]} g_{B!}
    \operatorname{in}_{\varphi!} \pr_{k,l}^* E[l] } }
\ar[r]<1ex>^(.63){g_E'} &
{ \displaystyle{ \coprod_{[k]} \pr_k^* E[k] }. } \cr
}$$
The map $g_{E,\theta}'$ is defined to be the following composite map:
$$g_{B!} \operatorname{in}_{\theta!} \pr_{m,n}^*E[n] \stackrel{\sim}{\ra}
\operatorname{in}_{[n]!} (\id \times \theta_{\Delta})_!
\pr_{m,n}^*E[n] \ra \operatorname{in}_{[n]!} \pr_n^* E[n].$$
The first map is the unique natural isomorphism
$g_{B!}\operatorname{in}_{\theta!} \stackrel{\sim}{\ra}
\operatorname{in}_{[n]!}(\id \times \theta_{\Delta})_!$ that exists because
$$g_B \circ \operatorname{in}_{\theta} = \operatorname{in}_{[n]} \circ
(\id \times \theta_{\Delta}) \colon B[n] \times \Delta^m \ra
\coprod_{[k]} B[k] \times \Delta^k.$$
The second map is induced from the map
$$(\id \times \theta_{\Delta})_! \pr_{m,n}^* E[n] \ra \pr_n^* E[n]$$
that is the adjoint of the unique natural isomorphism
$$\pr_{m,n}^* E[n] \stackrel{\sim}{\ra}
(\id \times \theta_{\Delta})^* \pr_n^* E[n]$$
which exists because
$$\pr_{m,n} = \pr_n \circ (\id \times \theta_{\Delta}) \colon B[n]
\times \Delta^m \ra B[n].$$  Then $g_E$ is defined to be the map
$\epsilon_{B!}g_E'$.  
The parametrized orthogonal spectrum $E$ over $B$ has the following mapping property, the proof of which follows directly from the definition of the parametrized orthogonal spectrum $E$ over $B.$
\begin{prop} Let $X$ be a parametrized orthogonal
spectrum over $B$. Then giving a map $\alpha \colon E \ra X$ of
parametrized orthogonal spectra over $B$ is equivalent to giving, for
every non-negative integer $k$, a map of parametrized orthogonal
spectra over $B[k] \times \Delta^k$
$$\alpha_k \colon \pr_k^*E[k] \ra
(\epsilon_B \circ \operatorname{in}_{[k]})^* X$$
such that, for every map $\theta \colon [m] \ra [n]$ in the simplicial
index category,  
$$(\theta_B \times \id)^* \alpha_m = (\id \times \theta_{\Delta})^*
\alpha_n \colon \pr_{m,n}^*E[n] \ra (\epsilon_B' \circ \operatorname{in}_{\theta})^* X.$$
\end{prop}
To understand the $\lambda$-th space $E_\lambda$ of the parametrized orthogonal spectrum $E$ over $B$ we note that limits, colimits, and the functors $f_!$, $f^\ast$, and $f_\ast$ are defined levelwise.  The space $E_\lambda$ is therefore given by the following coequalizer diagram of parametrized spaces over $B$:
$$\xymatrix{
{ \displaystyle{\epsilon_{B!}' (\coprod_{\theta \colon [m] \ra
      [n]}^{\phantom{\theta \colon [m] \ra [n]}} \pr_{m,n}^*E[n]_{\lambda}) } } \ar[r]<.7ex>^(.55){f_{E,\lambda}}
\ar[r]<-.7ex>_(.55){g_{E,\lambda}} &
{ \displaystyle{\epsilon_{B!} (\coprod_{ [k] }^{\phantom{ [k] }}
    \pr_k^*E[k]_{\lambda} ) } }
\ar[r]^(.7){\epsilon_{E,\lambda}} &
{ E_{\lambda}. } \cr
}$$
We recall the space $\theta_B^*E[m]_{\lambda}$ is defined by a
pullback diagram. Hence, there is a canonical map of spaces 
$\operatorname{pr} \colon \theta_B^*E[m]_{\lambda} \ra E[m]_\lambda$, and we define
$$\theta_{E,\lambda}^{\#} = \operatorname{pr} \circ\, \theta_{E,\lambda} \colon
E[n]_{\lambda} \ra E[m]_{\lambda}.$$
The map $\theta_{E,\lambda}^{\#}$ is a map of spaces and the following
diagrams commute:
$$\xymatrix{
{ E[n]_{\lambda} } \ar[r]^{\theta_{E,\lambda}^{\#}} \ar[d]^{p} &
{ E[m]_{\lambda} } \ar[d]^{p} &
{ E[n]_{\lambda} } \ar[r]^{\theta_{E,\lambda}^{\#}} &
{ E[m]_{\lambda} } \cr
{ B[n] } \ar[r]^{\theta_B} &
{ B[m] } &
{ B[n] } \ar[r]^{\theta_B} \ar[u]^{s} &
{ B[m]. } \ar[u]^{s} \cr 
}$$
Let $E_{\lambda}'$ be the space defined by the following
coequalizer diagram:
$$\xymatrix{
{ \displaystyle{\coprod_{\theta \colon [m] \to
      [n]}^{\phantom{\theta \colon [m] \to [n]}} \pr_{m,n}^*E[n]_{\lambda} } } \ar[r]<.7ex>^(.55){f_{E,\lambda}^{\#}}
\ar[r]<-.7ex>_(.55){g_{E,\lambda}^{\#}} &
{ \displaystyle{ \coprod_{ [k] }^{\phantom{ [k] }}
   \pr_k^* E[k]_{\lambda} } }
\ar[r]^(.67){\epsilon_{E,\lambda}^{\#}} &
{ E_{\lambda}'. } \cr
}$$
Then there is a canonical map $E_{\lambda}' \to E_{\lambda}$ and we
claim that this map is a homeomorphism.  Indeed, this is a special case
of the following more general statement:
\begin{prop}   Let $X[-]$ and $B[-]$ be two
diagrams of spaces indexed by a small category $I$, and let $X =
\colim_I X[-]$ and $B = \colim_I B[-]$. Let $p[-] \colon X[-] \ra
B[-]$ and $s[-] \colon B[-] \ra X[-]$ be natural transformations such
that $p[-] \circ s[-]$ is the identity natural transformation of
$B[-]$. Then each canonical map $\iota_{B,\alpha}\colon B[\alpha] \ra B$ gives rise to a parametrized space $(\iota_{B,\alpha})_! X[\alpha]$ over $B$ and the induced map $$\colim_I (\iota_{B,\alpha})_! X[\alpha] \ra X$$ is an isomorphism of parametrized spaces over $B$, where the induced maps $p \colon X \ra B$ and $s \colon B \ra
X$ provide the projection and section maps for the parametrized space $X$ over $B$. \end{prop}
\begin{proof}  We have commutative
diagrams
$$\xymatrix{
{ X[\alpha] } \ar[r]^{\iota_{X,\alpha}} \ar[d]^{p[\alpha]} &
{ X } \ar[d]^{p} &
{ X[\alpha] } \ar[r]^{\iota_{X,\alpha}} &
{ X } \cr
{ B[\alpha] } \ar[r]^{\iota_{B,\alpha}} &
{ B } &
{ B[\alpha] } \ar[r]^{\iota_{B,\alpha}} \ar[u]_{s[\alpha]} &
{ B. } \ar[u]_{s} \cr
}$$
It follows that we obtain a parametrized space $(\iota_{B,\alpha})_!X[\alpha]$
over $B$ together with a map $\tilde{\iota}_{X,\alpha} \colon
(\iota_{B,\alpha})_!X[\alpha] \to X$ of parametrized spaces over $B$.
The parametrized spaces $(\iota_{B,\alpha})_!X[\alpha]$ over $B$ form
an $I$-diagram of parametrized spaces over $B$, and the maps
$\tilde{\iota}_{X,\alpha}$ give rise to the following map of
parametrized spaces over $B$:
$$\tilde{\iota} \colon \colim_I (\iota_{B,\alpha})_!X[\alpha] \to X.$$
Now the general statement is that this map is an isomorphism of
parametrized spaces over $B$. Indeed, the canonical maps
$$\varphi_{\alpha} \colon X[\alpha] \to (\iota_{B,\alpha})_!X[\alpha]
\to \colim_I (\iota_{B,\beta})_!X[\beta]$$
give rise to a map
$$\varphi \colon X = \colim_I X[\alpha] \to \colim_I
(\iota_{B,\alpha})_!X[\alpha]$$
which is the inverse of the map $\tilde{\iota}$. 
\end{proof}
\begin{cor} Let $\omega\colon B \ra \ast$ be the map from $B$ to the one-point space and let $\omega_k\colon B[k] \ra \ast$ be the map from $B[k]$ to the one-point space.  The Thom spectrum $\omega_!E$ of the parametrized orthogonal spectrum $E$
over $B$ is canonically isomorphic to the realization of the
simplicial orthogonal spectrum given by the Thom spectra of the
parametrized orthogonal spectra $E[k]$ over $B[k]$, $$\omega_!E \stackrel{\sim}{\ra} \big| [k] \mapsto \omega_{k!} E[k]\big|$$ with simplicial
structure maps induced from the maps $\theta_{E,\lambda}^{\#}$.
\end{cor}

\section{Definition and Structure of the Parametrized Orthogonal Spectrum $E(A,G)$}\label{defofE(A,G)}

In general, for any symmetric ring spectrum $E$ as defined in \cite{hoveyshipleysmith}, we can define the topological Hochschild homology spectrum $\THH(E)$ following the approach of B\"{o}kstedt \cite{bokstedt}.  Further details of this construction may be found in \cite[\S 1--\S2]{h3}.
We now recall the definition of $\THH(E).$  Initially we define an index category $I$ by declaring the objects $ob(I)$ to be the class of all finite sets, ${\underline i} = \{1,2, \cdots , i \}, i \ge 1$ and ${\underline 0}= \emptyset $.  We then declare our morphisms to be {\it all} injective maps.  We note that every morphism is the composite of
the standard inclusion $m \mapsto m$ and an automorphism of the target set, albeit non-uniquely.  We now have, for each symmetric ring spectrum $E$ and pointed space $X$, a functor $G_{k}(E;X):I^{k+1}\ra \mathcal{T}_{\ast}$.  This is defined 
on objects as
\begin{center}
$G_{k}(E;X)({\underline i_{0}}, \dots , {\underline i_{k}}) = F(S^{i_{0}} \wedge \cdots \wedge S^{i_{k}}, E_{i_{0}} \wedge \cdots \wedge E_{i_{k}} \wedge X)$
\end{center} the space of based maps in $\mathcal{T}_{\ast}$ with the compact-open topology.  We define the functor on morphisms as follows.  Let $\iota \colon {\underline i_{r}}\ra {\underline i'_{r}}$ be the standard inclusion and write $i'_{r} = i_r + j_r.$  Then the natural transformation $G_{k}(E;X)({\underline i_{0}}, \dots , \iota, \dots,  {\underline i_{k}})$ takes the map $$S^{i_{0}}\wedge \cdots \wedge S^{i_r}\wedge \cdots \wedge S^{i_{k}} \stackrel{f}{\ra} E_{i_{0}} \wedge \cdots \wedge E_{i_r} \wedge \cdots \wedge E_{i_{k}}\wedge X$$ to the composition
\begin{align*} 
S^{i_{0}}\wedge \cdots \wedge S^{i'_r} \wedge \cdots \wedge \cdots \wedge S^{i_{k}} & \xrightarrow[\phantom{\mu_{i_r,j_r}}]{\sim} S^{i_{0}}\wedge \cdots \wedge S^{i_{r}}  \wedge \cdots \wedge S^{i_{k}}\wedge S^{j_r}\\
{}&\xrightarrow[\phantom{\mu_{i_r,j_r}}]{f\wedge \eta_E} E_{i_{0}}\wedge \cdots \wedge E_{i_{r}} \wedge \cdots \wedge E_{i_{k}}\wedge X \wedge E_{j_r}\\
{}&\xrightarrow[\phantom{\mu_{i_r,j_r}}]{\sim} E_{i_{0}}\wedge \cdots \wedge E_{i_{r}} \wedge E_{j_r} \wedge \cdots \wedge E_{i_{k}}\wedge X\\
{}&\xrightarrow{\mu_{i_r,j_r}} E_{i_{0}}\wedge \cdots \wedge E_{i'_{r}} \wedge \cdots \wedge E_{i_{k}}\wedge X 
\end{align*} where the first and third maps are the canonical isomorphisms.  The symmetric group $\Sigma_{i'_r}$ acts on $S^{i'_r}$ and $E_{i'_r}$ and by conjugation on $G_{k}(E;X)({\underline i_{0}}, \dots , {\underline i'_r}, \dots,  {\underline i_{k}})$.  This completes the definition of the functor $G_k$.
We define a cyclic pointed space with $k$-simplices $$\THH(E;X)[k] = \operatornamewithlimits{hocolim}_{I^{k+1}} G_{k}(E;X).$$  We define the pointed space $\THH(E;X)$ as the realization of the above cyclic pointed space
$$\THH(E;X) = \big| [k]\mapsto \THH(E;X)[k] \big|.$$
Let $\lambda$ be a finite dimensional inner product space and let $S^{\lambda}$ be the one point compactification
of $\lambda.$  Then we define the orthogonal spectrum $\THH(E)$ with $\lambda$-th space $$\THH(E)_{\lambda} = \THH(E; S^{\lambda}) = \big| [k] \mapsto \THH(E;S^{\lambda})[k]\big|$$ where $$\THH(E;S^{\lambda})[k] = \operatornamewithlimits{hocolim}_{I^{k+1}}G_{k}(E;S^{\lambda}).$$  For fixed $k$ and varying $\lambda$, each $\THH(E)[k]$ forms an orthogonal spectrum and $\THH(E)$ is the geometric realization of the resulting cyclic orthogonal spectrum.  
To define the face operators $d_{r}:\THH(E;S^{\lambda})[k]\ra \THH(E;S^{\lambda})[k-1]$, let $$I\times I \stackrel{\sqcup}{\ra} I$$ be the concatenation functor sending $(\underline{i}, \underline{i'})$ to the set $\underline{i}\sqcup \underline{i'}.$  Let $\partial_{r}:I^{k+1} \ra I^{k}$ be the functor defined by 
$$\partial_{r}({\underline i_{0}}, \dots , {\underline i_{k}}) = \begin{cases}({\underline i_{0}}, \dots ,{\underline i_{r}} \sqcup {\underline i_{r+1}}, \dots , {\underline i_{k}}) & \text{if $0\leq r < k$} \\
({\underline i_{k}}\sqcup {\underline i_{0}}, {\underline i_{1}}, \dots , {\underline i_{k-1}}) & \text{if $r=k.$}\end{cases} $$
Similarly for the degeneracies and cyclic operator, let ${\it s}_{r}: I^{k+1} \ra I^{k+2}$ be the functor defined on objects by $${\it s}_{r}({\underline i_{0}}, \dots , {\underline i_{k}}) = ({\underline i_{0}}, \dots ,  {\underline i_{r}}, {\underline 0}, {\underline i_{r+1}}, \dots , {\underline i_{k}})$$ for $0\leq r \leq k,$ and ${\it t}_{k}:I^{k}\ra I^{k}$ by
$${\it t}_{k} ({\underline i_{0}}, \dots , {\underline i_{k}}) = ({\underline i_{k}}, {\underline i_{0}}, \dots , {\underline i_{k-1}}).$$  We then define natural transformations 
\begin{align*} \delta_{r}& : G_{k}(E;S^{\lambda}) \ra G_{k-1}(E;S^{\lambda})\circ \partial_{r} \\
\sigma_{r} & :G_{k}(E;S^{\lambda})\ra G_{k+1}(E;S^{\lambda})\circ {\it s}_{r}\\
\tau_{k} & :G_{k}(E;S^{\lambda})\ra G_{k}(E;S^{\lambda})\circ {\it t}_{k} \end{align*}
as follows.  The natural transformation $\delta_{r}$ takes the map $f \in G_{k}(E;S^{\lambda}) ({\underline i_{0}}, \dots , {\underline i_{k}})$ given by $$S^{i_{0}}\wedge \cdots \wedge S^{i_{k}} \stackrel{f}{\ra} E_{i_{0}} \wedge \cdots \wedge E_{i_{k}}$$ to the map $\delta_{r}(f)\in G_{k-1}(E;S^{\lambda})(\partial_{r}({\underline i_{0}}, \dots , {\underline i_{k}}))$ given by the composition
\begin{align*} 
S^{i_{0}}\wedge \cdots \wedge S^{i_{r}+i_{r+1}}\wedge \cdots \wedge S^{i_{k}} & \stackrel{\sim}{\ra} S^{i_{0}}\wedge \cdots \wedge S^{i_{r}} \wedge S^{i_{r+1}} \wedge \cdots \wedge S^{i_{k}}\\
{}&\stackrel{f}{\ra} E_{i_{0}}\wedge \cdots \wedge E_{i_{r}} \wedge E_{i_{r+1}} \wedge \cdots \wedge E_{i_{k}}\wedge S^{\lambda} \\
{}&\stackrel{\mu}{\ra} E_{i_{0}}\wedge \cdots \wedge E_{i_{r}+i_{r}+1} \wedge \cdots \wedge E_{i_{k}}\wedge S^{\lambda}
\end{align*} if $0\leq r < k$ and 
\begin{align*}
S^{i_{k}+i_{0}}\wedge S^{i_{1}}\wedge \cdots \wedge S^{i_{k-1}} & \stackrel{\sim}{\ra} S^{i_{0}}\wedge \cdots \wedge S^{i_{k-1}}\wedge S^{i_{k}}\\
{}&\stackrel{f}{\ra} E_{i_{0}}\wedge \cdots \wedge E_{i_{k-1}} \wedge E_{i_{k}}\wedge S^{\lambda} \\
{}&\stackrel{\sim}{\ra} E_{i_{k}}\wedge E_{i_{0}} \cdots \wedge E_{i_{k-1}} \wedge S^{\lambda}\\
{}&\stackrel{\mu}{\ra} E_{i_{k}+i_{0}}\wedge E_{i_{1}} \cdots \wedge E_{i_{k-1}}\wedge S^{\lambda}
\end{align*} if $r=k.$  The natural transformations $\sigma_{r}$ and $\tau_{k}$ are defined similarly.  The face map $d_{r}:\THH(E;S^{\lambda})[k]\ra \THH(E;S^{\lambda})[k-1]$ is then defined as the 
composite $$\operatornamewithlimits{hocolim}_{I^{k+1}} G_{k}(E;S^{\lambda}) \stackrel{\delta_{r}}{\ra} \operatornamewithlimits{hocolim}_{I^{k+1}} G_{k-1}(E;S^{\lambda})\circ \partial_{r} \stackrel{(\partial_{r})_{\ast}}{\ra} \operatornamewithlimits{hocolim}_{I^{k}} G_{k-1}(E;S^{\lambda}).$$

For any ring $A$, we have an associated symmetric ring spectrum $\tilde{A}$, the {\it Eilenberg--MacLane spectrum}, with level $n$ obtained as the realization of the following simplicial set:
\begin{displaymath} 
\tilde{A}_{n}= {\big |} [k] \mapsto A\{{S}^{n}[k]\}/A\{s_{0}[k]\} {\big |},
\end{displaymath}
where we put the usual simplicial structure on the sphere,
\begin{center}
$S^{n}[-] = (\Delta^{1}[-] / \partial \Delta^{1}[-]) \wedge \cdots \wedge (\Delta^{1}[-]/ \partial \Delta^{1}[-])$ 
\end{center} with $n$ smash factors and basepoint $s_{0}[-] \in S^{n}[-]$.  We then define the {\it topological Hochschild spectrum of the ring} $A$ to be the topological Hochschild spectrum of the Eilenberg--MacLane spectrum $\tilde{A}$ associated to the
ring $A$, and simply write $\THH(A)$ for this spectrum.

\begin{defn}\label{actiondef}  Let $E$ be a connective symmetric ring spectrum and let $G$ be a discrete group.  A left action of $G$ on the spectrum $E$ is a continuous map
$$\alpha : G_{+} \wedge E_{n}  \ra E_{n}$$ such that the following diagrams commute
\[\xymatrix{G_{+} \wedge E_{m} \wedge E_{n}\ar[r]^{id \wedge \mu}\ar[d]^{\Delta \wedge id\wedge id} & G_{+} \wedge E_{m+n} \ar[d]^{\alpha}\\
G_{+} \wedge G_{+} \wedge E_{m}\wedge E_{n}\ar[d]^{\sim} & E_{m+n}\\
G_{+}\wedge E_{m} \wedge G_{+} \wedge E_n \ar[r]^(.6){\alpha \wedge \alpha} & E_{m}\wedge E_{n}\ar[u]_{\mu}},\] 
and
\[\xymatrix{G_{+}\wedge S^{m}\ar[d]^{pr}\ar[r]^(.47){id \wedge \eta}& G_{+}\wedge E_{m}\ar[d]^{\alpha}\\
S^{m}\ar[r]^{\eta}& E_{m}.}\]
\end{defn}  We note in the definition that $E$ is not necesarily the Eilenberg--MacLane spectrum of a ring and that we do not require the ring spectrum to be connected ($0$-connected) but only {\it connective} ($-1$-connected).  Also, the commutativaty of the these two diagrams implies the diagram 
\[\xymatrix{G_{+} \wedge E_{n}\wedge S^{m}  \ar[r]^(.55){\alpha \wedge id}\ar[d]^{id\wedge \sigma_{m,n}}& E_{n}\wedge S^{m}\ar[d]^{\sigma_{m,n}} \\
G_{+} \wedge E_{m+n} \ar[r]^(.55){\alpha}& E_{m+n}} \] commutes.  Here $\sigma_{m,n}$ is the structure map in the symmetric ring spectrum $E.$  We now define the {\it twisted group ring spectrum} to be the symmetric ring spectrum with level $n$ space given by $$(E^{\tau}[G])_{n} = E_{n}\wedge G_{+}.$$  We define the multiplication $(E^{\tau}[G])_{m}\wedge (E^{\tau}[G])_{n} \ra (E^{\tau}[G])_{m+n}$ as the composition 
\begin{align*}
E_{m}\wedge G_{+}\wedge E_{n} \wedge G_{+}&\xrightarrow{id \wedge \Delta \wedge id \wedge id}E_{m}\wedge G_{+}\wedge G_{+} \wedge  E_{n}\wedge G_{+}\\
{}&\xrightarrow[\phantom{id \wedge \Delta \wedge id \wedge id}]{id \wedge id \wedge \alpha \wedge id}E_{m}\wedge G_{+}\wedge E_{n} \wedge G_{+}\\
{}&\xrightarrow[\phantom{id \wedge \Delta \wedge id \wedge id}]{id \wedge tw \wedge id}E_{m}\wedge E_{n}\wedge G_{+} \wedge G_{+}\\
{}&\xrightarrow[\phantom{id \wedge \Delta \wedge id \wedge id}]{\mu_{E}\wedge \mu_{G}}E_{m+n}\wedge G_{+}
\end{align*}where $\mu_{E}$ and $\mu_{G}$ are the multiplication maps for $E$ and $G$ respectively.  We define the unit map $\eta: S^{m}\ra(E^{\tau}[G])_{m}$ to be the composite
$$S^{m}\xrightarrow[\phantom{\eta_{A}\wedge 1_{G}}]{\sim} S^{m}\wedge S^{0}\xrightarrow{\eta_{E}\wedge 1_{G}} E_{m}\wedge G_{+}$$ where $\eta_{E}:S^{m} \ra E_{m}$ is the unit map for the ring spectrum $E$ and $1_{G}:S^0\ra G_{+}$ is the constant map to the identity of the group $G.$

We relate this definition of a twisted group ring spectrum to the usual Eilenberg--MacLane spectrum of a twisted group ring by the following proposition.
\begin{prop}\label{twistedringprop}
Let $A$ be a ring, let $G$ be a discrete group, and let $\tilde{A}$ be the Eilenberg--Lane spectrum associated to $A$.  Then there exists a canonical weak equivalence of ring spectra $$\tilde{A}^{\tau}[G] \ra \widetilde{A^{\tau}[G]}.$$ 
\end{prop}
\begin{proof}
The canonical map is given by the composition $$(\tilde{A}^{\tau}[G])_{n} = \tilde{A}_{n} \wedge G_{+} \xrightarrow{\tilde{\phi}\wedge id} (\widetilde{A^{\tau}[G]})_{n}\wedge G_{+} \xrightarrow[\phantom{\tilde{phi}\wedge id}]{\tilde{r}} (\widetilde{A^{\tau}[G]})_{n}$$ where $\tilde{\phi}$ is induced from the ring homomorphism $\phi: A\ra A^{\tau}[G]$ defined by $\phi(a)= a \cdot 1$ and the map $\tilde{r}(-\wedge h)$ is induced from the ring homomorphism $r_{h}: A^{\tau}[G]\ra A^{\tau}[G]$ defined by $r_{h}(a\cdot g) = a \cdot gh.$ 
This composition induces isomorphisms of spectrum homotopy groups
and is thus a weak equivalence.
\end{proof}

Let $(g_{0}, \dots, g_{k})$ be a tuple
of elements of $G.$  We define a map $$(g_{0}, \dots, g_{k})_{\ast} : \THH(A)[k] \ra \THH(A)[k].$$  Let $f\in G_{k}(A)({\underline i_{0}}, \dots , {\underline i_{k}})$ be given by $$S^{i_{0}}\wedge \cdots \wedge S^{i_{k}} \stackrel{f}{\ra} \tilde{A}_{i_{0}} \wedge \cdots \wedge \tilde{A}_{i_{k}}.$$ We define a natural transformation $\gamma$ by sending $f$ to the composite
\begin{align*} 
S^{i_{0}}\wedge \cdots \wedge S^{i_{k}} & \stackrel{f}{\ra} \tilde{A}_{i_{0}}\wedge \cdots \wedge \tilde{A}_{i_{k}} \\
{}&\stackrel{\iota_{g}}{\ra} G_{+}\wedge \tilde{A}_{i_{0}} \wedge \cdots \wedge G_{+}\wedge \tilde{A}_{i_{0}}\\
{}& \ra \tilde{A}_{i_{0}}\wedge \cdots \wedge \tilde{A}_{i_{k}} 
\end{align*} where $\iota_{g}(a_{i_{0}} \wedge \cdots \wedge a_{i_{k}}) =  g_{0}\wedge a_{i_{0}} \wedge \cdots g_{k}\wedge a_{i_{k}},$ and the last map is $\alpha_{i_{0}}\wedge\cdots\wedge \alpha_{i_{k}}.$  Then $$(g_{0}, \dots, g_{k})_{\ast} : \THH(A)[k] \ra \THH(A)[k]$$ is the induced map
$$\operatornamewithlimits{hocolim}_{I^{k+1}}G_{k}(A)\stackrel{\gamma_{\ast}}{\ra}\operatornamewithlimits{hocolim}_{I^{k+1}}G_{k}(A).$$  

The topological Hochschild spectrum is defined as the geometric realization $$\THH(A)_{\lambda} = \big| [k] \mapsto \THH(A;S^{\lambda})[k] \big|$$ of the cyclic orthogonal spectrum $\THH(A,S^{\lambda})[-]$ with cyclic operators $d_{r}, s_{r},$ and $t_{k}$ as defined above.  We now define the cyclic orthogonal spectrum $\THH^{g}(A)$ that depends on a choice of element $g\in G.$  For $g\in G$ let $\THH^{g}(A)[-]$ be the cyclic orthogonal spectrum with $k$-simplices $$\THH^{g}(A)[k]_{\lambda} 
= \THH(A;S^{\lambda})[k]$$ and cyclic structure maps 
$$d_{r}^{g} = \begin{cases} d_{0}\circ (1,g, 1, \dots 1)_{\ast}&\text{if $r=0,$}\\
d_{r} &\text{if $0< r \leq k,$}
\end{cases}$$ $s_{r}^{g}= s_{r},$ and $t_{k}^{g}= t_{k}.$  The geometric realization is the orthogonal spectrum with $\lambda$-th space $$\THH^{g}(A)_{\lambda} = \big| [k] \mapsto \THH^{g}(A)[k]_{\lambda}\big|.$$  

The cyclic bar construction of a group $G$ is the cyclic set $N^{\cy}(G)[-]$ with $k$-simplices $$N^{\cy}(G)[k] = \underbrace{G\times \cdots 
\times G}_{k+1}$$ and face and degeneracy maps given by 
\begin{align*}d_{i}(g_{0}, \dots, g_{k}) &= \begin{cases} (g_{0}, \dots, g_{i}g_{i+1}, \dots, g_{k}) & \text{if $0 \leq i < k$}\\
(g_{k}g_{0}, g_{1}, \dots, g_{k-1})&\text{if $i=k$}
\end{cases}\\
s_{i}(g_{0}, \dots, g_{k}) &= (g_{0}, \dots, g_{i}, 1, g_{i+1}, \dots, g_{k}) \text{   for $0 \leq i \leq k$}.
\end{align*}  We also have the cyclic operator, $t_{k}$, defined by
$$t_{k}(g_{0}, \dots, g_{k}) = (g_{k}, g_{0}, \dots, g_{k-1}).$$

For each non-negative integer $k,$ let $E(A,G)[k]$ be the parametrized orthogonal spectrum over $N^{\cy}(G)[k]$ given by $$E(A,G;S^{\lambda})[k] = {\omega_k}^{\ast}(\THH(A;S^{\lambda})[k]) = \THH(A;S^{\lambda})[k]\times N^{\cy}(G)[k]$$ where ${\omega_k}$ is the unique map from $N^{\cy}(G)$ to the one-point space.  To define the structure maps of the parametrized orthogonal spectrum $E(A,G),$ we first recall the spectrum structure maps of the orthogonal spectrum $\THH(A)[k]$.  The space $\THH(A;S^\lambda)[k]$ is obtained as the homotopy colimit of spaces of the form $F(X, Y \wedge S^\lambda),$ and the spectrum structure map 
$$\sigma_{\lambda,\lambda'} \colon \THH(A;S^{\lambda})[k] \wedge
S^{\lambda'} \ra \THH(A;S^{\lambda \oplus \lambda'})[k]$$
is then obtained from the canonical map
$$F(X,Y \wedge S^{\lambda}) \wedge S^{\lambda'} \ra F(X,Y \wedge
S^{\lambda} \wedge S^{\lambda'})$$
and various canonical isomorphisms. It is clear that this makes
$\THH(A)[k]$ an orthogonal spectrum.  We now define the {\it twisted} structure maps
$$\theta_E^{\tau} \colon E(A,G)[n] \ra
\theta_B^*E(A,G)[m]$$
of parametrized orthogonal spectra over $B[n] = N^{\cy}(G)[n]$.  First we note that we have corresponding untwisted structure maps
$$\theta_E \colon E(A,G)[n] \ra
\theta_B^*E(A,G)[m]$$
defined by $\theta_E = \omega_n^*\theta_{\THH(A)}$. We define
$$\begin{aligned}
s_{r,E}^{\tau} = s_{r,E} & \colon E(A,G)[k] \ra
s_{r,B}^* E(A,G)[k+1], \hskip 4mm 0 \leqslant r \leqslant k, \cr 
t_{k,E}^{\tau} = t_{k,E} & \colon E(A,G)[k] \ra E(A,G)[k] \cr
\end{aligned}$$
to be the untwisted maps, and define
$$d_{r,E}^{\tau} = d_{r,E} \circ \varphi_r \colon E(A,G)[k] \ra
d_{r,B}^* E(A,G)[k-1], \hskip 4mm 0 \leqslant r \leqslant k,$$
to be the composition of the untwisted map and the automorphism
$$\varphi_r \colon E(A,G)[k] \ra E(A,G)[k]$$
of parametrized orthogonal spectra over $N^{\cy}(G)[k]$ defined by
$$\varphi_r(f,(g_0,\dots,g_k)) =
((t_k^{r+1}(g_r,1,\dots,1))_*(f),(g_0,\dots,g_k)).$$
We verify that the map $\varphi_r$ is a map of
parametrized orthogonal spectra over $B[k] = N^{\cy}(G)[k]$; that is, we
check that the following diagram of parametrized spaces over
$B[k]$ commutes:
$$\xymatrix{
{ E(A,G;S^{\lambda})[k] \wedge_{B[k]} S_{B[k]}^{\lambda'} }
\ar[rr]^{\varphi_{r,\lambda} \wedge \id}
\ar[d]^{\sigma_{\lambda,\lambda'}} &&
{ E(A,G;S^{\lambda})[k] \wedge_{B[k]} S_{B[k]}^{\lambda'} }
\ar[d]^{\sigma_{\lambda,\lambda'}} \cr
{ E(A,G;S^{\lambda \oplus \lambda'})[k] }
\ar[rr]^{\varphi_{r,\lambda \oplus \lambda'}} &&
{ E(A,G;S^{\lambda \oplus \lambda'})[k]. } \cr
}$$
Since the commutativity of the above diagram is established by verifying the commutativity of the maps at the point-set level, we can check the diagram one fiber at a time. The induced diagram of fibers over
$(g_0,\dots,g_k) \in B[k]$ takes the form
$$\xymatrix{
{ \THH(A;S^{\lambda})[k] \wedge S^{\lambda'} }
\ar[rrr]^{(t_k^{r+1}(g_r,1,\dots,1))_* \wedge \id}
\ar[d]^{\sigma_{\lambda,\lambda'}} &&&
{ \THH(A;S^{\lambda})[k] \wedge S^{\lambda'} }
\ar[d]^{\sigma_{\lambda,\lambda'}} \cr
{ \THH(A;S^{\lambda \oplus \lambda'})[k] }
\ar[rrr]^{(t_k^{r+1}(g_r,1,\dots,1))_*} &&&
{ \THH(A;S^{\lambda \oplus \lambda'})[k].} \cr
}$$
It commutes since the canonical map
$$F(X,Y \wedge S^{\lambda}) \wedge S^{\lambda'} \ra F(X,Y \wedge
S^{\lambda} \wedge S^{\lambda'})$$
is natural in the variables $X$ and $Y$. This completes the definition
of the twisted structure maps. 

We next show that given composable maps $\theta \colon [m] \ra [n]$ and $\theta'
\colon [n] \ra [p]$ in the simplicial index category, the following diagram of parametrized orthogonal
spectra over $B[p]$ commutes:
$$\xymatrix{
{ E(A,G)[p] }
\ar[rr]^(.4){\theta_E^{\tau} \circ \theta_E'{}^{\hskip-1mm \tau}} 
\ar[d]^{\theta_E'{}^{\hskip-1mm\tau}} &&
{ (\theta_B \circ \theta_B')^* E(A,G)[m] } \ar[d]^{\sim} \cr
{ \theta_B'{}^* E(A,G)[n] } \ar[rr]^(.45){\theta_B'{}^*\theta_E^{\tau}} &&
{ \theta_B'{}^*\theta_B^* E(A,G)[m]. } \cr
}$$
Again, this can be easily be checked on fibers. Hence, we obtain a
parametrized orthogonal spectrum $E(A,G)$ over $N^{\cy}(G)$ with
$\lambda$th space
$$E(A,G)_\lambda = E(A,G;S^{\lambda}) = \big| [k] \mapsto E(A,G;S^{\lambda})[k] \big|$$
where the simplicial structure maps in the simplicial space on the
right-hand side are the twisted maps $\theta_E^{\tau \#}$.  We now present the proof of Theorem ~\ref{secondmainresult}.
\begin{proof}[Proof of Theorem ~\ref{secondmainresult}] Let $\omega$ be the unique map from $B = N^{\cy}(G)$ to the one-point
space. We wish to construct a map of parametrized orthogonal spectra
over $B$
$$\tilde{\Psi} \colon E(A,G) \ra \omega^*\THH(A^{\tau}[G])$$
and show that the adjoint map
$$\Psi \colon \omega_! E(A,G) \ra \THH(A^{\tau}[G])$$
is a stable equivalence of orthogonal spectra. These maps exist for
every symmetric ring spectrum $R$ with a $G$-action in the sense of
Def.~\ref {actiondef} and with the symmetric ring spectrum $R^{\tau}[G]$ as defined
in the paragraph immediately following Def.~\ref{actiondef}. We shall work in this generality.  As we noted above, the orthogonal spectrum $\omega_! E(R,G)$ is the realization of the simplicial orthogonal spectrum with $k$-th term $${\omega_k}_! E(R,G)[k] = \THH(R)[k] \wedge N^{\cy}(G)[k]_+$$ and with simplicial structure maps given by the compositions $${\theta}_{\omega_!E}^\tau \colon {\omega_n}_!E(R,G)[n]\xrightarrow[\phantom{{\omega_m}_!\epsilon \omega^*_m}]{{\omega_n}_!\theta_E^\tau}{\omega_n}_! \theta_B^*E(R,G)[m]\xrightarrow{{\omega_m}_!\epsilon \omega^*_m} {\omega_m}_!E(R,G)[m].$$  Here, we recall, ${\omega_n}_! = {\omega_m}_! {\theta_B}_!$.  The map $\epsilon: {\theta_B}_!\theta^*_B \ra \id$ is the counit of the adjunction $({\theta_B}_!, \theta^*_B)$ and is given by the map $$\epsilon = \id \wedge \theta_B \colon \THH(R)[m] \wedge N^{\cy}(G)[n]_+ \ra \THH(R)[m] \wedge N^{\cy}(G)[m]_+.$$ 
The map $\Psi$ is defined to be the map of realizations obtained from a map of simplicial orthogonal spectra $$\Psi_k \colon \THH(R)[k] \wedge N^{\cy}(G)[k]_+ \ra \THH(R^{\tau}[G])[k]$$ that we define below.  The definition of this map is given in the proof of~\cite[Thm.~7.1]{hm}. It is also shown there that the map is a stable
equivalence of orthogonal spectra (provided that the symmetric ring
spectrum $R$ converges; this is the case for $R = \tilde{A}$). Hence,
it suffices to show that the maps $\Psi_k$ are compatible with the
simplicial structure maps. We first recall that the map
$$\Psi_{k,\lambda} \colon \THH(R;S^{\lambda})[k] \wedge
N^{\cy}(G)[k]_+ \ra \THH(R^{\tau}[G];S^{\lambda})[k].$$
is the map of homotopy colimits over $I^{k+1}$ obtained from the
composite map
$$\begin{aligned}
{} & F(S^{i_0} \wedge \dots \wedge S^{i_k},R_{i_0} \wedge \dots
R_{i_k} \wedge S^{\lambda}) \wedge G_+ \wedge \dots \wedge G_+ \cr
{} & \ra F(S^{i_0} \wedge \dots \wedge S^{i_k},R_{i_0} \wedge \dots
\wedge R_{i_k} \wedge G_+ \wedge \dots \wedge G_+ \wedge S^{\lambda}) \cr
{} & \ra F(S^{i_0} \wedge \dots \wedge S^{i_k},R_{i_0} \wedge
G_+ \wedge \dots \wedge R_{i_k} \wedge G_+ \wedge S^{\lambda}) \cr
\end{aligned}$$
where the first map is the same canonical map that was used to define
the spectrum structure maps, and where the second map is induced from
the permutation
$$R_{i_0} \wedge \dots \wedge R_{i_k} \wedge G_+ \wedge \dots \wedge 
G_+ \ra R_{i_0} \wedge G_+ \wedge \dots \wedge R_{i_k} \wedge G_+$$
that maps $(r_0,\dots,r_k,g_0,\dots,g_k)$ to $(r_0,g_0,\dots,r_k,g_k)$.  We show the following diagram commutes:
\[\xymatrix{{\omega_k}_! E(R,G)\ar[r]^(.475){\Psi_k}\ar[d]^{d^\tau_{r,\omega_! E}}& \THH(R^\tau[G])[k]\ar[d]^{d_r}\\
{\omega_{(k-1)}}_! E(R,G)[k-1]\ar[r]^{\Psi_{k-1}}&\THH(R^\tau[G])[k-1].
}\]  
Let $\partial_r \colon I^{k+1} \ra I^k$ be the functor defined in the beginning of the section and let $G_k(R;X)$ be the functor from $I^{k+1}$ to the category of pointed spaces also defined at beginning of the section.  Let $\delta_r \colon G_k(R;X) \ra G_{k-1}(R;X) \circ \partial_r$ be the natural transformation used to define the face map of the cyclic pointed space $\THH(R;X)[-]$, again, defined at the beginning of the section.  Then the right-hand vertical map in the above diagram is given by the composite map 
\begin{align*} \hocolim_{I^{k+1}} G_k(R^\tau[G]; S^\lambda) & \xrightarrow[\phantom{\partial_{r\ast}}]{\delta_{r\ast}} \hocolim_{I^{k+1}} G_{k-1}(R^\tau[G];S^\lambda )\circ \partial_r\\
{}& \xrightarrow[\phantom{\delta_{r\ast}}]{\partial_{r\ast}} \hocolim_{I^k} G_{k-1} (R^\tau[G];S^\lambda). \end{align*}  The left-hand vertical map in the diagram above is also a composition, given by
\begin{align*} \hocolim_{I^{k+1}} G_k(R; S^\lambda)\wedge N^{\cy}(G)[k]_+ & \xrightarrow[\phantom{\partial_{r\ast}}]{\delta^\tau_{r\ast}} \hocolim_{I^{k+1}} G_{k-1}(R;S^\lambda )\wedge N^{\cy}(G)[k-1]_+\circ \partial_r\\
{}& \xrightarrow[\phantom{\delta^\tau_{r\ast}}]{\partial_{r\ast}} \hocolim_{I^k} G_{k-1} (R;S^\lambda) \wedge N^{\cy}(G)[k-1]_+ \end{align*} where the natural transformation $$\delta^\tau_r \colon G_k(R;X) \wedge N^{\cy}(G)[k]_+\ra (G_{k-1}(R;X) \wedge N^{\cy}(G)[k-1]_+)\circ \partial_r$$ is defined by $$\delta^\tau_r (f, (g_0, \dots, g_k)) = (\delta_r((t^{r+1}_k(g_r, 1, \dots, 1))_\ast (f)), d_r (g_0, \dots, g_k)).$$  Hence, it suffices to show that the diagram of natural transformations 
\[\xymatrix{G_k(R;S^\lambda)\wedge N^{\cy}(G)[k]_+ \ar[rr]^(.55){\Psi_{k,\lambda}}\ar[d]^{\delta^\tau_r}&&G_k(R^\tau[G];S^\lambda)\ar[d]^{\delta_r}\\
G_{k-1}(R;S^\lambda) \wedge N^{\cy}(G)[k-1]_+ \circ \partial_r \ar[rr]^(.575){\Psi_{k-1, \lambda} \circ \partial_r} && G_{k-1}(R^\tau[G];S^\lambda) \circ \partial_r
}\] commutes.  But this follows immediately from the definitions of the natural transformations involved and from the naturality of the canonical map $$F(X, Y) \wedge Z \ra F(X, Y \wedge Z).$$  We therefore have the desired map $\Psi$ of orthogonal spectra and its adjoint $\tilde{\Psi}$ of parametrized orthogonal spectra over $N^{\cy}(G)$.  As we previously stated, it is proved in~\cite[Thm.~7.1]{hm} (see also ~\cite[Prop.~4.1]{h3}) that the maps of orthogonal spectra $\Psi_k$ are stable equivalences, provided that the symmetric ring spectrum $R$ converges.  We wish to also conclude that the induced map of realizations $\Psi$ is a stable equivalence.  It is proved in ~\cite{may} that this holds, provided the simpicial spaces $\omega_! E(R,G)[-]$ and $\THH(R^\tau[G])[-]$ are {\it proper} in the sense of ~\cite[Def. 11.2]{may}.  This, in turn, is the case, if the unit maps $\eta_i \colon S^i \ra R_i$ are Hurewicz cofibrations.  If $R= \tilde{A}$, then both properties hold, and hence, the map $\Psi$ is a stable equivalence of orthogonal spectra.  
 \end{proof}

\section{The Fiber Bundle $E(A,G)_\lambda$ over $N^{\cy}(G)$}\label{THHcyc}

The cyclic set $EG[-]$ is defined by $EG[k]=\operatorname{Map}([k],G).$  The face and degeneracy operators $d_{i}:EG[k]\ra EG[k-1]$ and $s_{i}:EG[k]\ra EG[k+1]$, for $0\leq i \leq k,$ are given by 
\begin{align*}
d_{i}(g_{0}, \dots, g_{k}) &= (g_{0}, \dots, \hat{g_{i}}, \dots, g_{k})\\
s_{i}(g_{0}, \dots, g_{k}) &= (g_{0}, \dots, g_{i}, g_{i}, \dots, g_{k}),
\end{align*}where the hat symbol indicates that the $i$-th term is omitted.  The cyclic operator is defined $$t_{k}(g_{0},\dots, g_{k}) = (g_{k}, g_{0},\dots, g_{k-1}).$$  We follow \cite{schlichtkrull}, and let $G^{\operatorname{ad}}$ denote the set $G$ with the group $G$ acting from the left by conjugation.  We note that an element $g\in G$ determines an isomorphism of sets between $G/C_{G}(g)$ and the 
conjugacy class $\langle g\rangle$ of the element $g$ given by mapping the class $aC_{G}(g)$ to $g^{a} = aga^{-1}.$  Here $C_{G}(g)$ denotes the centralizer of $g.$  Therefore, $$G^{\operatorname{ad}} = \coprod\limits_{\langle g\rangle} G/C_{G}(g) \cdot g.$$  Then  
the map $$\phi:EG[-]\times_{G} G^{\operatorname{ad}} \stackrel{\sim}{\ra} N^{\cy}(G)[-]$$ defined on level $k$ by
$$\phi([(g_{0}, \dots, g_{k}); g]) = (g_{k}gg_{0}^{-1}, g_{0}g_{1}^{-1}, g_{0}g_{1}^{-1}, \dots, g_{k-1}g_{k}^{-1})$$ is an isomorphism.  The inverse sends a $k$-simplex $(g_{0}, \dots, g_{k})$ in $N^{\cy}(G)[k]$ to
the class $[(g_{1}\cdots g_{k}, g_{2}\cdots g_{k}, \dots,  g_{k}, 1 ); g_{0}g_{1}\cdots g_{k}].$
The adjoint of the composition $$\mathbb{T}\times N^{\cy}(G)\stackrel{\mu}{\ra}N^{\cy}(G)\stackrel{\pi}{\ra}BG$$  is a map to the free-loop space of the classifying space of the group, $$N^{\cy}(G)\ra \Lambda BG$$ and this map is a weak equivalence \cite[~Prop 2.6]{bokstedthsiangmadsen}.  
The map $$\pi: N^{\cy}(G)\ra N(G) = BG$$ in the composition is given by the projection $$(g_0,\dots, g_{k}) \mapsto (g_{1}, \dots, g_{k}).$$  Also, the set of connected components of $N^{\cy}(G)$ is in one-to-one correspondence with the set of conjugacy classes of elements in the group $G$.
Given the cyclic space $E(A,G)[-]_{\lambda}$ defined in \S\ref{defofE(A,G)}, we now describe the geometric fiber over the connected component of $E(A,G)_{\lambda}$ corresponding to the conjugacy class of the element $g \in G$.
To understand the fiber over the point given by the $0$-simplex $g \in N^{\cy}(G)[0],$ we evaluate the pullback
\[\xymatrix{E(A,G)_{\lambda}^{g} \ar[d] \ar[r]& E(A,G)_{\lambda} \ar[d]\\
\Delta^{0} \ar[r]^(.4){g} & N^{\cy}(G)}\]
where the right-hand vertical map is the projection of the cyclic space $E(A,G)[-]_{\lambda}$ onto the cyclic bar construction and the bottom horizantal map
takes the unique non-degenerate $0$-simplex to $g$.  Recall that geometric realization preserves finite limits \cite{gabrielzisman} in the sense that the 
canonical map $$\big| [k] \mapsto \operatornamewithlimits{lim}_{\alpha} X_{\alpha} [-] \big| \ra \operatornamewithlimits{lim}_{\alpha} \big| [k] \mapsto X_{\alpha}\big|$$ is a homeomorphism, provided that the index category for the limit 
system is finite.  Hence the geometric fiber $E(A,G)^{g}$ is the space obtained as the realization of the following pullback diagram of simplicial spaces  
\[\xymatrix{E(A,G)^{g}[-]_{\lambda} \ar[d] \ar[r]& E(A,G)[-]_{\lambda} \ar[d]\\
\Delta^{0}[-] \ar[r]^(.4){g} & N^{\cy}(G)[-]. }\]
Here $\Delta^{0}[-]$ is the standard $0$-simplex $$\Delta^{0} =\operatorname{Hom}_{\Delta}(-, [0])$$ and the bottom map takes the identity, $[0] \stackrel{{\bf 1}}{\ra} [0]$ to the $0$-simplex $g \in N^{\cy}(G)[0]$.  Explicitly, this map takes the map $\theta:[k] \ra [0]$ to $\theta^{\ast}g \in N^{\cy}(G)[k],$
and since for each $k$, there is only one such map $\theta$, $\theta^{\ast}g = (g,1, \dots, 1)$ ($k+1$ factors).  Thus the fiber $E(A,G)^{g}[-]_{\lambda}$ is given simplical degree-wise by $$E(A,G)^{g}[k]_{\lambda}= \THH(A;S^{\lambda})[k] \times \{(g,1, \dots, 1)\} \subset E(A,G)_{\lambda}[k].$$  The cyclic structure of the fiber is that of $E(A,G)[-]_{\lambda}$ restricted at each $k$ to the subset $\{(g,1,\dots, 1)\} \subset N^{\cy}(G)[k].$  Thus as a simplicial set the fiber is canonically isomorphic to the simplicial set $\THH^{g}(A)[-]$ defined in \S\ref{defofE(A,G)}.
The connected component of the $0$-simplex $g$ is obtained as the realization of the cyclic spectrum $E(A,G)\big{|}_{\langle g \rangle}[-]_{\lambda}$ given on level $k$ by\\ $$E(A,G)\big|_{\langle g \rangle}[k]_{\lambda} = \THH(A;S^{\lambda})[k] \times \big{\{}(g_{0}, g_{1}, \dots, g_{k}) : g_0g_1\cdots g_k \in \langle g \rangle \big{\}},$$ a subset of $E(A,G)[k]_{\lambda}.$  Indeed, a path from $g$ to the zeroth vertex of the $k$-simplex $(g_{0}, \dots, g_{k})$ is given by the $1$-simplex
$(gh^{-1},h)\in N^{\cy}(G)[1].$  
\begin{lem}\label{nutsnbolts} There exists an equivalence of orthogonal spectra $$\tilde{\phi}_{g}:EG\times_{C_{G}(g)} \THH^{g}(A)\stackrel{\sim}{\ra} E(A,G)\big{|}_{\langle g\rangle},$$ between the Borel construction and the spectrum corresponding to the connected component indexed by $\langle g \rangle.$  The equivalence depends on the choice of representative $g\in \langle g \rangle.$
\end{lem}
\begin{proof}
A stronger statement holds.  Namely, we have a degree-wise isomorphism  $$\tilde{\phi}_{g}: EG[-]\times_{C_{G}(g)} \THH^{g}(A)[-] \ra E(A,G)\big{|}_{\langle g\rangle}[-]_{\lambda}$$ between the Borel construction $EG[-]\times_{C_{G}(g)} \THH^{g}(A)[-]_{\lambda}$ and the cyclic space $E(A,G)\big{|}_{\langle g\rangle}[-]_{\lambda}.$  The isomorphism is defined by
$$\tilde{\phi}_{g}([(g_{0}, \dots, g_{k}); [f]]) = ((g_{k},g_{0}, \dots, g_{k-1})_{\ast}(f); \phi(g_{0}, \dots, g_{k})),$$ where $\phi(g_{0}, \dots, g_{k}) = (g_{k}gg_{0}^{-1}, g_{0}g_{1}^{-1}, \dots, g_{k-1}g_{k}^{-1})$ as above and $(f)$ is the class in the homotopy colimit represented by the map  $$S^{i_{0}}\wedge \cdots S^{i_{k}} \stackrel{f}{\ra} \tilde{A}_{i_{0}} \wedge \cdots \wedge \tilde{A}_{i_{k}}.$$
It is readily verified that the map $\tilde{\phi}_{g}$ is well-defined and using Def.~\ref{actiondef} it is also easily checked that the map commutes with the face operators $d_{r}\times d^{g}_{r}$ and $d^{\tau}_{r,E}.$  The other cyclic operators of $EG[-]\times_{C_{G}(g)} \THH^{g}(A)[-]$ are the standard product operators of $EG[-]$ and $\THH(A)[-]_{\lambda}.$  These operators do not involve the element $g\in G,$ and their commutavity with their corresponding maps $s^{\tau}_{r,E}$ and $t^{\tau}_{k,E}$ is even more easily verified.
\end{proof}

We state Lemma~\ref{nutsnbolts} globally as follows.  Let $X[-]_{\lambda}$ be the cyclic space with $k$-simplices $X[k]_{\lambda} = G^{\operatorname{ad}} \times \THH^{g}(A)[k]_{\lambda}$ and cyclic operators those of $\THH^{g}(A)[-]_{\lambda}.$  As a set, $X[k]_{\lambda}$ is the disjoint union of all fibers $$X[k]_{\lambda}=\coprod_{\langle g \rangle}\THH^{g}(A)[k]_{\lambda}.$$  Now given the Borel construction $EG[-] \times_{G} X[-]_{\lambda}$ with the usual product
simplicial structure, we give a cyclic structure by defining the cyclic operator 
$$t(g_{0}, \dots, g_{k}; g, [f]) = (g_{k}g, g_{0}, g_{1}, \dots, g_{k-1}; g, \tau_{k}\circ (1, g, \cdots, 1)_{\ast}[f]),$$ where $\tau_{k}$ is the natural transformation from \S\ref{defofE(A,G)}.  For a $k$-simplex in $EG[-] \times_{G} X[-]_{\lambda}$, we define the map $$\tilde{\phi}: EG[-] \times_{G} X[-]_{\lambda} \ra E(A,G)[-]_{\lambda}$$ via
$$\tilde{\phi}(\bar{g} ; g, [f]) = ((g_{k}, g_{0}, g_{1}, \dots, g_{k-1})_{\ast}[f]; \phi_{g}),$$ where $\bar{g} = (g_{0}, \dots, g_{k})$ and $\phi_{g} = (g_{k}gg_{0}^{-1}, g_{0}g_{1}^{-1}, g_{0}g_{1}^{-1}, \dots, g_{k-1}g_{k}^{-1}).$  We note that the restriction of this map to the connected component corresponding to $\langle g \rangle$ and a choosen representative $g\in G$ is 
the map $\tilde{\phi}_{g}$ of Lemma~\ref{nutsnbolts}.  This map is an isomorphism and we see at once we have the following proposition.
\begin{prop} There exist canonical homeomorphisms of spaces $\phi$ and $\tilde{\phi}$ such that the following diagram commutes:
\[\xymatrix{ EG\times_{G} X_{\lambda} \ar[r]^(.525){\tilde{\phi}} \ar[d]& E(A,G)_{\lambda}\ar[d] \\
EG \times_{G} G^{\operatorname{ad}} \ar[r]^(.55){\phi} & N^{\cy}(G).}\] 
\end{prop}
\begin{proof}  The cyclic isomorphism $\phi$ is covered by the cyclic isomorphism $\tilde{\phi};$ that is we have a diagram of cyclic spaces:
\[\xymatrix{ EG[-]\times_{G} X[-]_{\lambda} \ar[r]^(.55){\tilde{\phi}} \ar[d]& E(A,G)[-]_{\lambda}\ar[d] \\
EG[-] \times_{G} G^{\operatorname{ad}} \ar[r]^(.55){\phi} & N^{\cy}(G)[-]}\] where the two horizontal maps are cyclic isomorphisms.  After taking geometric realization, the result follows.
\end{proof}
\begin{cor}\label{maincor} There are canonical homeomorphisms of spaces:
\[\xymatrix{ \coprod\limits^{\phantom{\langle g \rangle}}_{\langle g\rangle}EG\times_{G} \THH^{g}(A)_{\lambda} \ar[r]^(.65){\sim} \ar[d]& E(A,G)_{\lambda}\ar[d] \\
EG  \times_{G} \coprod\limits^{\phantom{\langle g \rangle}}_{\langle g \rangle} G/C_{G}(g) \cdot g \ar[r]^(.65){\sim} & N^{\cy}(G)}\] where the two horizontal maps are the homeomorphisms corresponding to the cyclic isomorphisms.
\end{cor}  
\begin{rem}  We note that the space $E(A,G)_{\lambda}$ is actually a fiber bundle over $N^{\cy}(G)$.  In particular, the fibers over two points in the same connected component are homeomorphic.
\end{rem}  We now prove Theorem~\ref{mainresult}.
\begin{proof}[Proof of Theorem ~\ref{mainresult}]  Applying the functor $f_{!}$ to the parametrized spaces over $N^{\cy}(G),$ $$\coprod\limits_{\langle g \rangle} EG\times_{G} E(A,G)_{\lambda}^{g}\ra EG \times_{G} \coprod\limits_{\langle g \rangle} G/C_{G}(g) \stackrel{\sim}{\ra}N^{\cy}(G)$$ and 
$$E(A,G)_{\lambda}\ra N^{\cy}(G)$$ gives a map of spaces $$\bigvee\limits_{\langle g \rangle} EG\wedge_{C_{G}(G)} E(A,G)_{\lambda}^{g}\stackrel{\sim}{\ra}f_{!}E(A,G)_{\lambda}\ra\THH(A^{\tau}[G])_{\lambda}.$$ The first map is an isomorphism by Corollary~\ref{maincor}.  As $\lambda$ varies, the second map is an equivalence by Theorem~\ref{secondmainresult}.  Hence for varying $\lambda,$ we obtain a stable equivalence of orthogonal spectra.
\end{proof}

\providecommand{\bysame}{\leavevmode\hbox to3em{\hrulefill}\thinspace}
\providecommand{\MR}{\relax\ifhmode\unskip\space\fi MR }
\providecommand{\MRhref}[2]{%
  \href{http://www.ams.org/mathscinet-getitem?mr=#1}{#2}
}
\providecommand{\href}[2]{#2}

\bibliographystyle{amsplain}
\end{document}